\documentclass[11pt]{article}%
\usepackage{amsmath}
\usepackage{amsfonts}
\usepackage{amssymb}
\usepackage{graphicx}
\usepackage{a4wide}%
\providecommand{\U}[1]{\protect\rule{.1in}{.1in}}
\newtheorem{theorem}{Theorem}

\newtheorem{corollary}[theorem]{Corollary}

\newtheorem{example}[theorem]{Example}

\newtheorem{lemma}[theorem]{Lemma}

\newtheorem{proposition}[theorem]{Proposition}

\begin{document}

\title{On the lower semi\-continuity and subdifferentiability of the value function
for conic linear programming problems}
\author{Constantin Z\u{a}linescu\thanks{Octav Mayer Institute of Mathematics, Ia\c{s}i
Branch of Romanian Academy, Ia\c{s}i, Romania, and University ``Alexandru Ioan
Cuza" Ia\c{s}i, Romania; email: \texttt{zalinesc@uaic.ro}.}}
\date{}
\maketitle

\begin{abstract}
Lemma 1 from the paper [N.E. Gretsky, J.M. Ostroy, W.R. Zame,
Subdifferentiability and the duality gap, Positivity 6: 261--274, 2002]
asserts that the value function $v$ of an infinite dimensional linear
programming problem in standard form is lower semicontinuous whenever $v$ is
proper and the involved spaces are normed vector spaces. In this note one
shows that this statement is false even in finite-dimensional spaces, one
provides an example of linear programming problem in Hilbert spaces whose
(proper) value function is not lower semi\-continuous (hence it is not
sub\-differentiable) at any point in its domain, one shows that the
restriction of the value function to its domain in Kretschmer's gap example is
not bounded on any neighborhood of any point of the domain, and discuss other
assertions done in the same paper.

\end{abstract}

\section{Introduction}

The following conical linear programming problem

\medskip(P) \ minimize $\ c^{\ast}(x)$ \ s.t. \ $x\in P,\ Ax-b\in Q$,

\medskip\noindent and its dual

\medskip(D) \ maximize $\ y^{\ast}(b)$ \ s.t. $\ y^{\ast}\in Q^{+}$, $c^{\ast
}-A^{\ast}y^{\ast}\in P^{+}$,

\medskip\noindent are studied in \cite{GrOsZa02}, where $X$, $Y$ are Hausdorff
locally convex spaces, $X^{\ast}$ and $Y^{\ast}$ are their topological dual
spaces, $A:X\rightarrow Y$ is a continuous linear operator, $A^{\ast}:Y^{\ast
}\rightarrow X^{\ast}$ is the adjoint of $A$, $P\subset X$ and $Q\subset Y$
are convex cones, $P^{+}\subset X^{\ast}$ and $Q^{+}\subset Y^{\ast}$ are the
positive dual cones of $P$ and $Q$, $b\in Y$ and $c\in X^{\ast}$.

\smallskip The main results from \cite{GrOsZa02} are: Theorem 1 which states
that the value function $v$ associated to (P) is subdifferentiable at $b$ if
and only if (D) has optimal solutions and there is no duality gap, and its use
for proving the Duffin--Karlovitz no-gap theorem; Lemma 1 which states that
$v$ is lower semicontinuous whenever it is proper and the involved spaces are
normed vector spaces; the modification of Kretschmer's gap example in order to
get a convex function which is subdifferentiable at a point but is not
continuous there; Proposition 2 which provides sufficient conditions, adequate
for the assignment model, to ensure that $v$ is Lipschitz on $Q$.

Unfortunately, \cite[Lem.\ 1]{GrOsZa02} is not true even in finite dimensional
spaces, which makes its use to be not adequate in the proof of the
Duffin--Karlovitz no-gap theorem, while the proof of \cite[Prop.\ 2]{GrOsZa02}
needs, in our opinion, serious clarifications; moreover, there are also other
inaccuracies in the paper. Having in view the remark that \textquotedblleft
This paper should be on the reading list of any advanced mathematical
economics course which has a focus on extremal methods in infinite-dimensional
spaces" from the review MR1932651 (2003i:90111) of \cite{GrOsZa02} in
Mathematical Reviews, we consider that there is a strong motivation for an
attentive reading of this paper.

\smallskip

The paper is organized as follows. Having in view that the value function
associated to problem (P) is positively homogeneous and subadditive, in
Section 2 we underline some specific properties of such functions, pointing
out the differences between the (lower-, upper-, local Lipschitz) continuity
of such functions and their restrictions to the domain. In Section 3 we
essentially discuss the proof of the Duffin--Karlovitz no-gap theorem, while
in Section 4 we provide two counter-examples to \cite[Lem.\ 1]{GrOsZa02}, the
first in finite-dimensional spaces and the second in Hilbert spaces. Section 5
is dedicated to Kretschmer's gap example, while in Section 6 we comment the
proof of \cite[Prop.\ 2]{GrOsZa02}.

\smallskip Below, we introduce the basic notations and some preliminary
results used in the paper.

\smallskip Throughout this note, the considered spaces are real Hausdorff
locally convex spaces (H.l.c.s.\ for short) if not mentioned explicitly
otherwise. Having $X$ an H.l.c.s., $X^{\ast}$ is its topological dual endowed
with its weakly-star topology $w^{\ast}:=\sigma(X^{\ast},X)$. The value
$x^{\ast}(x)$ of $x^{\ast}\in X^{\ast}$ at $x\in X$ is denoted by
$\left\langle x,x^{\ast}\right\rangle $. It is well known that $(X^{\ast
},w^{\ast})^{\ast}$ can be identified with $X$, what we do in the sequel.
Having $(\emptyset\neq)$ $K\subset X$ a convex cone (that is, $x+x^{\prime}\in
K$ and $tx\in K$ for all $x,x^{\prime}\in K$ and $t\in\mathbb{R}%
_{+}:=[0,\infty\lbrack$), we set $x\leq_{K}x^{\prime}$ (equivalently
$x^{\prime}\geq_{K}x$) for $x,x^{\prime}\in X$ with $x^{\prime}-x\in K$;
clearly $\leq_{K}$ is a \emph{preorder} on $X$, that is, $\leq_{K}$ is
reflexive and transitive. For $E\subset X$, one denotes by
$\operatorname*{span}E$, $\operatorname*{aff}E$, $\operatorname*{icr}E$,
$\operatorname*{cor}E$, $\operatorname*{int}E$ and $\operatorname*{cl}E$ the
\emph{linear hull}, the \emph{affine hull}, the \emph{intrinsic core}, the
\emph{core}, the \emph{interior} and the \emph{closure} of $E$, respectively.
For $\emptyset\neq A\subset X$ (and similarly for $\emptyset\neq B\subset
X^{\ast}$) we set $A^{+}:=\{x^{\ast}\in X^{\ast}\mid\forall x\in
A:\left\langle x,x^{\ast}\right\rangle \geq0\}$ for the \emph{positive dual
cone} of $A$; it is well known that $A^{+}$ $(\subset X^{\ast})$ is a
$w^{\ast}$-closed convex cone and $(K^{+})^{+}=\operatorname*{cl}K$ whenever
$K$ is a convex cone.

\smallskip Having a function $f:X\rightarrow\overline{\mathbb{R}}%
:=\mathbb{R}\cup\{-\infty,\infty\}$, its \emph{domain} is the set
$\operatorname*{dom}f:=\{x\in X\mid f(x)<\infty\}$; $f$ is \emph{proper} if
$\operatorname*{dom}f\neq\emptyset$ and $f(x)\neq-\infty$ for all $x\in X$;
$f$ is \emph{convex} if its \emph{epigraph} $\operatorname*{epi}f:=\{(x,t)\in
X\times\mathbb{R}\mid f(x)\leq t\}$ is convex; $f$ is \emph{positively
homogeneous} if $f(tx)=tf(x)$ for all $t\in\mathbb{P}:=\mathbb{R}_{+}%
\setminus\{0\}$ and $x\in X$; $f$ is \emph{sub\-additive} if $f(x+x^{\prime
})\leq f(x)+f(x^{\prime})$ for all $x,x^{\prime}\in\operatorname*{dom}f$; $f$
is \emph{sublinear} if $f$ is positively homogeneous, sub\-additive and
$f(0)=0$; $f$ is \emph{lower semi\-continuous} (l.s.c.\ for short) \emph{at}
$x\in X$ if $\liminf_{x^{\prime}\rightarrow x}f(x^{\prime})\geq f(x)$, where
$\overline{\mathbb{R}}$ is endowed with its usual topology (for example, the
topology induced by the metric $d$ defined by $d(t,t^{\prime}):=\left\vert
\arctan t-\arctan t^{\prime}\right\vert $ with $\arctan(\pm\infty):=\pm\pi
/2$); $f$ is \emph{l.s.c.}\ if $f$ is l.s.c.\ at any $x\in X$; the
\emph{l.s.c.\ envelope} of $f$ is the function $\overline{f}:X\rightarrow
\overline{\mathbb{R}}$ such that $\operatorname*{epi}\overline{f}%
=\operatorname*{cl}(\operatorname*{epi}f),$ and so $\overline{f}$ is convex if
$f$ is so; $f$ is \emph{upper semi\-continuous} (u.s.c.\ for short) \emph{at}
$x\in X$ (resp.\ on $X$) if $-f$ is l.s.c.\ at $x\in X$ (resp.\ on $X$); the
\emph{sub\-differential} of $f$ at $x\in X$ with $f(x)\in\mathbb{R}$ is the
set
\[
\partial f(x):=\{x^{\ast}\in X^{\ast}\mid\forall x^{\prime}\in X:\left\langle
x^{\prime}-x,x^{\ast}\right\rangle \leq f(x^{\prime})-f(x)\}
\]
and $\partial f(x):=\emptyset$ if $f(x)\notin\mathbb{R}$; $f$ is
\emph{sub\-differentiable} at $x\in X$ if $\partial f(x)\neq\emptyset$. The
\emph{conjugate} of $f$ is the function $f^{\ast}:X^{\ast}\rightarrow
\overline{\mathbb{R}}$ defined by
\[
f^{\ast}(x^{\ast}):=\sup\{\left\langle x,x^{\ast}\right\rangle -f(x)\mid x\in
X\}=\sup\{\left\langle x,x^{\ast}\right\rangle -f(x)\mid x\in
\operatorname*{dom}f\}\quad(x^{\ast}\in X^{\ast}),
\]
where $\sup\emptyset:=-\infty$; clearly, $f^{\ast}$ is a $w^{\ast}%
$-l.s.c.\ convex function. Having $g:X^{\ast}\rightarrow\overline{\mathbb{R}}%
$, its conjugate $g^{\ast}:X\rightarrow\overline{\mathbb{R}}$ of is defined
similarly. Notice that $f^{\ast}=(\overline{f})^{\ast}$; moreover, for $x\in
X$ and $x^{\ast}\in X^{\ast}$ one has%
\begin{equation}
x^{\ast}\in\partial f(x)\Leftrightarrow\lbrack f(x)\in\mathbb{R}%
\ \wedge\ f(x)+f^{\ast}(x^{\ast})=\left\langle x,x^{\ast}\right\rangle
]\Rightarrow\overline{f}(x)=f(x)\in\mathbb{R}\Rightarrow\partial\overline
{f}(x)=\partial f(x). \label{r-goz2}%
\end{equation}

The \emph{indicator function} of $E\subset X$ is $\iota_{E}:X\rightarrow
\overline{\mathbb{R}}$ defined by $\iota_{E}(x):=0$ for $x\in E$ and
$\iota_{E}(x):=\infty$ for $x\in X\setminus E$; notice that $\iota_{E}$ is
l.s.c.\ iff $E$ is closed, and $\iota_{E}$ is convex iff $E$ is convex.

\section{Some properties of sublinear functions\label{sec2}}

Because the value function of a linear programming problem is positively
homogeneous and sub\-additive (hence sublinear when it vanishes at $0$), it is
useful to point out some specific properties of such functions.

Let $g:X\rightarrow\overline{\mathbb{R}}$ be positively homogeneous and
sub\-additive. First observe that
\[
\forall x,x^{\prime}\in\operatorname*{dom}g,\ \forall\lambda\in\mathbb{P}%
\ :\ x+x^{\prime}\in\operatorname*{dom}g\ \ \text{and}\ \ \lambda
x\in\operatorname*{dom}g;\label{r-goz13}%
\]
consequently, for $x,x^{\prime}\in\operatorname*{dom}g$ and $\lambda\in
{}]0,1[$ one gets
\[
g(\lambda x+(1-\lambda)x^{\prime})\leq g(\lambda x)+g((1-\lambda)x^{\prime
})=\lambda g(x)+(1-\lambda)g(x^{\prime}),
\]
and so $g$ is convex.

Because $g$ is positively homogeneous, one has $g(0)=g(\lambda0)=\lambda g(0)$
for $\lambda\in\mathbb{P}$, and so $g(0)\in\{-\infty,0,\infty\}$. Moreover, if
$g(x_{0})=-\infty$ for some $x_{0}\in X$, then $g(x+\lambda x_{0})\leq
g(x)+\lambda g(x_{0})=-\infty$, whence $g(x+\lambda x_{0})=-\infty$, for all
$\lambda\in\mathbb{P}$ and $x\in\operatorname*{dom}g$; consequently,
$g(x)=-\infty$ for all $x\in\operatorname*{dom}g$ if $g(0)=-\infty$.

Assume that $g(0)=\infty$; then $\widetilde{g}:X\rightarrow\overline
{\mathbb{R}}$ defined by $\widetilde{g}(x):=g(x)$ for $x\neq0$ and
$\widetilde{g}(0):=0$ is sublinear. Indeed, take $x^{\prime},x^{\prime\prime
}\in\operatorname*{dom}\widetilde{g}=\{0\}\cup\operatorname*{dom}g$. If
$x^{\prime},x^{\prime\prime}\in\operatorname*{dom}g$, then $x^{\prime
}+x^{\prime\prime}\in\operatorname*{dom}g$ (hence $x^{\prime}+x^{\prime\prime
}\neq0$), and so $\widetilde{g}(x^{\prime}+x^{\prime\prime})=g(x^{\prime
}+x^{\prime\prime})\leq g(x^{\prime})+g(x^{\prime\prime})=\widetilde
{g}(x^{\prime})+\widetilde{g}(x^{\prime\prime})$; if $x^{\prime}=0$ (and
similarly for $x^{\prime\prime}=0$), then $\widetilde{g}(x^{\prime}%
+x^{\prime\prime})=\widetilde{g}(x^{\prime\prime})=\widetilde{g}(x^{\prime
})+\widetilde{g}(x^{\prime\prime})$. Because $\widetilde{g}$ is clearly
positively homogeneous, $\widetilde{g}$ is sublinear.

In the rest of this section, $g:X\rightarrow\overline{\mathbb{R}}$ is a
sublinear function; hence $g(0)=0$. Using \cite[Th.\ 2.4.14]{Zal02} one
obtains that%
\begin{gather}
\lbrack\partial g(0)\neq\emptyset\Leftrightarrow g\ \text{is\ l.s.c.\ at
}0],\quad g^{\ast}=\iota_{\partial g(0)},\quad\operatorname*{dom}g^{\ast
}=\partial g(0),\label{r-goz18a}\\
\partial g(0)\neq\emptyset\Rightarrow\left[  \forall x\in X:\overline
{g}(x)=\sup\{\left\langle x,x^{\ast}\right\rangle \mid x^{\ast}\in\partial
g(0)\}=g^{\ast\ast}(x)\right]  ,\label{r-goz18b}\\
\forall x\in X:\partial g(x)=\{x^{\ast}\in\partial g(0)\mid\left\langle
x,x^{\ast}\right\rangle =g(x)\}. \label{r-goz18c}%
\end{gather}
Also note that [$g$ is u.s.c.\ at $0$]\ $\Leftrightarrow$
[$\operatorname*{dom}g=X$ and $g$ is continuous on $X$]; moreover, [$g$ is
l.s.c.\ at $0$]\ $\Leftrightarrow$ [$g|_{\operatorname*{dom}g}$ is l.s.c.\ at
$0$], where $\operatorname*{dom}g$ is endowed with its trace (induced) topology.

Assume now that $X$ is a normed vector space; then one also has:
\begin{align}
g\ \text{is l.s.c.\ at }0  &  \Leftrightarrow\lbrack\exists L>0,\ \forall x\in
X:g(x)\geq-L\left\Vert x\right\Vert ];\label{r-goz17a}\\
g\ \text{is u.s.c.\ at }0  &  \Leftrightarrow\lbrack\exists L>0,\ \forall x\in
X:g(x)\leq L\left\Vert x\right\Vert ]\label{r-goz17b}\\
&  \Leftrightarrow\operatorname*{dom}g=X\ \text{and\ }g\ \text{is
(}L\text{-)Lipschitz on }X. \label{r-goz17c}%
\end{align}
Indeed, the implications $\Leftarrow$ from (\ref{r-goz17a})--(\ref{r-goz17c})
are obvious because $g(0)=0$. Assume that $g$ is l.s.c.\ at $0$. Then there
exists $r>0$ such that $g(x)\geq-1$ for $x\in X$ with $\left\Vert x\right\Vert
\leq r$. Taking $x\in X\setminus\{0\}$ and $x^{\prime}:=\frac{r}{\left\Vert
x\right\Vert }x$, one has $\left\Vert x^{\prime}\right\Vert \leq r$, whence
$-1\leq g(x^{\prime})=\frac{r}{\left\Vert x\right\Vert }g(x)$, and so
$g(x)\geq-r^{-1}\left\Vert x\right\Vert $; hence the implication $\Rightarrow$
holds in (\ref{r-goz17a}). The proof of the implication $\Rightarrow$ from
(\ref{r-goz17b}) is similar. Assume now that $g(x)\leq L\left\Vert
x\right\Vert $ for $x\in X$; hence $\operatorname*{dom}g=X$. Taking $x\in X$,
one has $0=g(0)=g(x+(-x))\leq g(x)+g(-x)$, whence $g(x)\in\mathbb{R}$. Take
now $x,x^{\prime}\in X$; we (may) assume that $g(x)\geq g(x^{\prime})$. Then
\[
g(x)=g((x-x^{\prime})+x^{\prime})\leq g(x-x^{\prime})+g(x^{\prime})\leq
L\left\Vert x-x^{\prime}\right\Vert +g(x^{\prime}),
\]
whence $\left\vert g(x)-g(x^{\prime})\right\vert =g(x)-g(x^{\prime})\leq
L\left\Vert x-x^{\prime}\right\Vert $. Therefore, $g$ is $L$-Lipschitz .

In what concern the continuity and the upper semi\-continuity of
$g|_{\operatorname*{dom}g}$, in a similar way as for (\ref{r-goz17b}), one
obtains:%
\begin{gather*}
g|_{\operatorname*{dom}g}\ \ \text{is u.s.c.\ at }0\Leftrightarrow
\lbrack\exists L>0,\ \forall x\in\operatorname*{dom}g:g(x)\leq L\left\Vert
x\right\Vert ],\\
g|_{\operatorname*{dom}g}\ \ \text{is continuous\ at }0\Leftrightarrow
\lbrack\exists L>0,\ \forall x\in\operatorname*{dom}g:\left\vert
g(x)\right\vert \leq L\left\Vert x\right\Vert ].
\end{gather*}

The next example shows the big differences among the continuity properties of
the functions $g$ and $g|_{\operatorname*{dom}g}$.

\begin{example}
\label{ex6}Let $X$ be an infinite-dimensional normed vector space (n.v.s.\ for
short) and let $\varphi:X\rightarrow\mathbb{R}$ be a linear, not continuous
functional. Consider the proper sublinear functions $g_{1},g_{2}%
,g_{3}:X\rightarrow\overline{\mathbb{R}}$ defined by
\[
g_{1}:=\max\{0,\varphi\},\quad g_{2}(x):=\left\{
\begin{array}
[c]{ll}%
0 & \text{if }x\in\lbrack\varphi\leq0],\\
\infty & \text{if }x\in\lbrack\varphi>0],
\end{array}
\right.  \quad g_{3}(x):=\left\{
\begin{array}
[c]{ll}%
\varphi(x) & \text{if }x\in\lbrack\varphi\leq0],\\
\infty & \text{if }x\in\lbrack\varphi>0],
\end{array}
\right.
\]
where $[\varphi\leq0]:=\{x\in X\mid\varphi(x)\leq0\}$ and similarly for
$[\varphi>0]$ (and the like). TFAH:

\emph{(i)} $\operatorname*{dom}g_{1}=X$, $\partial g_{1}(x)=\{0\}$ for
$x\in\lbrack\varphi\leq0]$, $\partial g_{1}(x)=\emptyset$ for $x\in
\lbrack\varphi>0]$, $(g_{1})^{\ast}=\iota_{\{0\}}$, $\overline{g_{1}}=0$, and
so $g_{1}$ is l.s.c.\ at $x$ iff $x\in\lbrack\varphi\leq0]$; moreover, $g_{1}$
is not u.s.c.\ at each $x\in X$.

\emph{(ii)} $\operatorname*{dom}g_{2}=[\varphi\leq0]$, $\partial
g_{2}(x)=\{0\}$ for $x\in\lbrack\varphi\leq0]$, $\partial g_{2}(x)=\emptyset$
for $x\in\lbrack\varphi>0]$, $(g_{2})^{\ast}=\iota_{\{0\}}$, $\overline{g_{2}%
}=0$, and so $g_{2}$ is l.s.c.\ at $x$ iff $x\in\lbrack\varphi\leq0]$;
moreover, $g_{2}$ is u.s.c.\ at $x$ iff $x\in\lbrack\varphi>0]$ and
$g_{2}|_{\operatorname*{dom}g_{2}}$ is Lipschitz.

\emph{(iii)} $\operatorname*{dom}g_{3}=[\varphi\leq0]$, $\partial
g_{3}(x)=\emptyset$ for $x\in X$, $(g_{3})^{\ast}=\infty$, $\overline{g_{3}%
}=-\infty$, and so $g_{3}$ is not l.s.c.\ at each $x\in X$; moreover, $g_{3}$
is u.s.c.\ at $x$ iff $x\in\lbrack\varphi>0]$, $g_{3}|_{\operatorname*{dom}%
g_{3}}$ is u.s.c.\ at $x$ iff $x\in\lbrack\varphi=0]$, and $g_{3}%
|_{\operatorname*{dom}g_{3}}$ is not l.s.c.\ at each $x\in\operatorname*{dom}%
g_{3}$.
\end{example}

Proof. (i) It is clear that $\operatorname*{dom}g_{1}=X$ and $0\in\partial
g_{1}(0)$. Consider $x^{\ast}\in\partial g_{1}(0)$; then, obviously,
$\left\langle x,x^{\ast}\right\rangle \leq g_{1}(x)=0$ for $x\in\lbrack
\varphi=0]$, and so $\left\langle x,x^{\ast}\right\rangle =0$ for $x\in
\lbrack\varphi=0]$. Using \cite[Lem.\ 3.9]{Rud91}, there exists $\lambda
\in\mathbb{R}$ such that $x^{\ast}=\lambda\varphi$, and so $\lambda=0$ because
$\varphi$ is not continuous. Hence $\partial g_{1}(0)=\{0\}$, whence,
$\partial g_{1}(x)=\{0\}$~$\Leftrightarrow$ $g_{1}(x)=0$~$\Leftrightarrow$
$\varphi(x)\leq0$, and $\partial g_{1}(x)=\emptyset$ when $\varphi(x)>0$ by
(\ref{r-goz18c}). Using (\ref{r-goz18a}) one gets $(g_{1})^{\ast}%
=\iota_{\{0\}}$, while using (\ref{r-goz18b}) one gets $\overline{g_{1}}=0$,
and so $g_{1}$ is l.s.c.\ at $x$ if and only if $\varphi(x)\leq0$. Assume that
$g_{1}$ is u.s.c.\ at $x\in X$; then $g_{1}$ is Lipschitz on
$\operatorname*{dom}g_{1}=X$ by (\ref{r-goz17c}), and so $g_{1}$ is l.s.c.\ on
$X$; this contradiction proves that $g_{1}$ is not u.s.c.\ at each $x\in X$.

(ii) It is obvious that $\operatorname*{dom}g_{2}=[\varphi\leq0]$ and
$0\in\partial g_{2}(0)$. Similar to the proof of (i) one gets that $\partial
g_{2}(x)=\{0\}$ if $x\in\lbrack\varphi\leq0]$ and $\partial g_{2}%
(x)=\emptyset$ otherwise, that $\overline{g_{2}}=0$, and so $g_{2}$ is
l.s.c.\ at $x$ if and only if $x\in\lbrack\varphi\leq0]$. Because
$g_{2}(x)=\infty\geq g_{2}(x^{\prime})$ for all $x\in\lbrack\varphi>0]$ and
$x^{\prime}\in X$, $g_{2}$ is obviously u.s.c.\ at $x\in\lbrack\varphi>0]$.
Assuming that $g_{2}$ is u.s.c.\ at some $x\in\operatorname*{dom}g_{2}$, one
obtains that $x\in\operatorname*{int}[\varphi\leq0]$, and so one gets the
contradiction that $\varphi$ is continuous; hence $g_{2}$ is u.s.c.\ at $x$
iff $x\in\lbrack\varphi>0]$. Because $g_{2}|_{\operatorname*{dom}g_{2}}=0$,
$g_{2}|_{\operatorname*{dom}g_{2}}$ is Lipschitz.

(iii) Clearly, $\operatorname*{dom}g_{3}=[\varphi\leq0]$. One has that
$x^{\ast}\in\partial g_{3}(0)$ iff $\left\langle x,x^{\ast}\right\rangle
\leq\varphi(x)$ for all $x\in\lbrack\varphi\leq0]$; using similar arguments to
those in the proof of (i), one gets $\partial g_{3}(0)=\emptyset$, and so
$\partial g_{3}(x)=\emptyset$ for every $x\in X$. Because $\partial
g_{3}(0)=\emptyset$, one has that $(g_{3})^{\ast}=\iota_{\emptyset}=\infty$,
and that $g_{3}$ is not l.s.c.\ at $0$, whence $\overline{g_{3}}(0)=-\infty$,
and so $\overline{g_{3}}(x)=-\infty$ for all $x\in\operatorname*{cl}%
(\operatorname*{dom}g_{3})=X$ because $[\varphi\leq0]$ is dense in $X$ as
$\varphi$ is not continuous. Hence $\overline{g_{3}}=-\infty$, and so $g_{3}$
is not l.s.c.\ at any $x\in X$. As in the proof of (ii), $g_{3}$ is u.s.c.\ at
$x$ iff $x\in\lbrack\varphi>0]$. Because $g_{3}|_{\operatorname*{dom}g_{3}%
}(x)\geq g_{3}|_{\operatorname*{dom}g_{3}}(x^{\prime})$ for all $x\in
\lbrack\varphi=0]$ and $x^{\prime}\in\operatorname*{dom}g_{3}$, $g_{3}%
|_{\operatorname*{dom}g_{3}}$ is u.s.c.\ at each $x\in\lbrack\varphi=0]$.

Because $\varphi$ is not continuous, $\sup\{\varphi(x)\mid x\in S_{X}%
\}=\infty$, where $S_{X}:=\{x\in X\mid\left\Vert x\right\Vert =1\}$, and so
there exists a sequence $(x_{n})_{n\geq1}\subset S_{X}$ such that
$\varphi(x_{n})\geq n$ for $n\in\mathbb{N}^{\ast}$. Consider $x\in
\lbrack\varphi<0]$ and take $x_{n}^{\prime}:=x-[\varphi(x)/\varphi
(x_{n})]x_{n}$ for $n\geq1$. Then $\varphi(x_{n}^{\prime})=0$ (whence
$x_{n}^{\prime}\in\operatorname*{dom}g_{3}$) for $n\in\mathbb{N}^{\ast}$,
$x_{n}^{\prime}\rightarrow x$, and $\limsup_{x^{\prime}\rightarrow x}%
\varphi(x^{\prime})\geq\limsup_{n\rightarrow\infty}\varphi(x_{n}^{\prime
})=0>\varphi(x)=g_{3}(x)$. Hence $g_{3}|_{\operatorname*{dom}g_{3}}$ is not
u.s.c.\ at $x$. Therefore, $g_{3}|_{\operatorname*{dom}g_{3}}$ is not
u.s.c.\ at each $x\in\lbrack\varphi<0]$, proving so that $g_{3}%
|_{\operatorname*{dom}g_{3}}$ is u.s.c.\ at $x\in\operatorname*{dom}g_{3}$ iff
$x\in\lbrack\varphi=0]$. Consider now $x\in\lbrack\varphi\leq0]$ and take
$x_{n}^{\prime\prime}:=x-[1/\varphi(x_{n})]x_{n}$ for $n\geq1$. Then
$\varphi(x_{n}^{\prime\prime})=\varphi(x)-1\leq0$ (whence $x_{n}^{\prime
\prime}\in\operatorname*{dom}g_{3}$) for $n\in\mathbb{N}^{\ast}$,
$x_{n}^{\prime\prime}\rightarrow x$, and $\liminf_{x^{\prime\prime}\rightarrow
x}\varphi(x^{\prime\prime})\leq\liminf_{n\rightarrow\infty}\varphi
(x_{n}^{\prime\prime})=\varphi(x)-1<\varphi(x)=g_{3}(x)$. Hence $g_{3}%
|_{\operatorname*{dom}g_{3}}$ is not l.s.c.\ at each $x\in\operatorname*{dom}%
g_{3}$. \hfill$\square$

\medskip

The next result seems to be quite relevant in the context of \cite{GrOsZa02}.

\begin{proposition}
\label{p-Z}Let $X$ be a n.v.s.\ and $g:X\rightarrow\overline{\mathbb{R}}$ be a
proper sublinear function. Assume that $x\in\operatorname*{dom}g$ and
$\delta,L>0$ are such that $\left\vert g(x^{\prime})-g(x^{\prime\prime
})\right\vert \leq L\left\Vert x^{\prime}-x^{\prime\prime}\right\Vert $ for
all $x^{\prime},x^{\prime\prime}\in B(x,\delta)\cap\operatorname*{dom}g$,
where $B(x,\delta):=\{x^{\prime}\in X\mid\left\Vert x^{\prime}-x\right\Vert
<\delta\}$. Then
\[
\forall\gamma\in\mathbb{P},\ \forall x^{\prime},x^{\prime\prime}\in B(\gamma
x,\gamma\delta)\cap\operatorname*{dom}g:\left\vert g(x^{\prime})-g(x^{\prime
\prime})\right\vert \leq L\left\Vert x^{\prime}-x^{\prime\prime}\right\Vert ;
\]
in particular, $g|_{\operatorname*{dom}g}$ is $L$-Lipschitz if $x=0$ (or, more
generally, $\left\Vert x\right\Vert <\delta$).
\end{proposition}

Proof. Take $x^{\prime},x^{\prime\prime}\in B(\gamma x,\gamma\delta
)\cap\operatorname*{dom}g$. Then $\gamma^{-1}x^{\prime},\gamma^{-1}%
x^{\prime\prime}\in B(x,\delta)\cap\operatorname*{dom}g$, and so
\[
\gamma^{-1}\left\vert g(x^{\prime})-g(x^{\prime\prime})\right\vert =\left\vert
g(\gamma^{-1}x^{\prime})-g(\gamma^{-1}x^{\prime\prime})\right\vert \leq
L\left\Vert \gamma^{-1}x^{\prime}-\gamma^{-1}x^{\prime\prime}\right\Vert
=L\gamma^{-1}\left\Vert x^{\prime}-x^{\prime\prime}\right\Vert ,
\]
whence $\left\vert g(x^{\prime})-g(x^{\prime\prime})\right\vert \leq
L\left\Vert x^{\prime}-x^{\prime\prime}\right\Vert $. \hfill$\square$

\medskip

Notice that $g|_{\operatorname*{dom}g}$ is locally Lipschitz on
$\operatorname*{icr}(\operatorname*{dom}g)$ whenever $g$ is a proper sublinear
function and $\dim X<\infty$ (or, more generally, $\dim X_{g}<\infty$, where
$X_{g}:=\operatorname*{span}(\operatorname*{dom}g)=\operatorname*{dom}%
g-\operatorname*{dom}g$); this is because for $\widetilde{g}:=g|_{X_{g}}%
\in\Lambda(X_{g})$, one has $\operatorname*{dom}\widetilde{g}%
=\operatorname*{dom}g$ and $\operatorname*{icr}(\operatorname*{dom}%
g)=\operatorname*{cor}(\operatorname*{dom}\widetilde{g})=\operatorname*{int}%
(\operatorname*{dom}\widetilde{g})$ and so the proper convex function
$\widetilde{g}$ is continuous (and so locally Lipschitz) on
$\operatorname*{int}(\operatorname*{dom}\widetilde{g})$.

\section{An alternative proof for Theorem 1 from \cite{GrOsZa02}}

If not mentioned explicitly otherwise, the problems (P) and (D), as well as
the corresponding data, are as in the Introduction. Moreover, the preorders
defined by $P$, $Q$, $P^{+}$ and $Q^{+}$ are simply denoted by $\leq$.

\medskip The value function associated to problem (P) is
\[
v:Y\rightarrow\overline{\mathbb{R}},\quad v(y):=\inf\{\left\langle x,c^{\ast
}\right\rangle \mid Ax\geq y,\ x\geq0\}=\inf\{\left\langle x,c^{\ast
}\right\rangle \mid x\in P,\ y\in Ax-Q\},
\]
where $\inf\emptyset:=\infty$. It is clear that $\operatorname*{dom}v=A(P)-Q$,
$v$ is positively homogeneous and convex, $v(0)\leq0$, and $v(y_{1})\leq
v(y_{2})$ whenever $y_{1}\leq y_{2}$; hence $v(0)\in\{0,-\infty\}$, and so
$v(y)=-\infty$ for $y\in\operatorname*{dom}v$ if $v(0)=-\infty$. Moreover, by
the definition of $v$, we get%
\begin{align*}
v^{\ast}(y^{\ast})  &  =\sup\big\{\left\langle y,y^{\ast}\right\rangle
+\sup\{\left\langle x,-c^{\ast}\right\rangle \mid x\in P,\ Ax-y=q\in Q\ \}\mid
y\in Y\big\}\\
&  =\sup\{\left\langle Ax-q,y^{\ast}\right\rangle +\left\langle x,-c^{\ast
}\right\rangle \mid x\in P,\ q\in Q\}\\
&  =\sup\{\left\langle x,A^{\ast}y^{\ast}-c^{\ast}\right\rangle +\left\langle
q,-y^{\ast}\right\rangle \mid x\in P,\ q\in Q\}\\
&  =\iota_{P^{+}}(c^{\ast}-A^{\ast}y^{\ast})+\iota_{Q^{+}}(y^{\ast})
\end{align*}
for every $y^{\ast}\in Y^{\ast}$; hence
\begin{equation}
\operatorname*{dom}v^{\ast}=\{y^{\ast}\in Q^{+}\mid c^{\ast}-A^{\ast}y^{\ast
}\in P^{+}\}=Q^{+}\cap(A^{\ast})^{-1}(c^{\ast}-P^{+}),\quad v^{\ast}%
=\iota_{\operatorname*{dom}v^{\ast}}, \label{r-goz1}%
\end{equation}
and so $\operatorname*{dom}v^{\ast}$ is the feasible set of problem (D). It
follows that
\[
v^{\ast\ast}(y)=\sup\left\{  \left\langle v,y^{\ast}\right\rangle -v^{\ast
}(y^{\ast})\mid y^{\ast}\in Y^{\ast}\right\}  =\sup\left\{  \left\langle
y,y^{\ast}\right\rangle \mid y^{\ast}\in Q^{+},\ c^{\ast}-A^{\ast}y^{\ast}\in
P^{+}\right\}
\]
for all $y\in Y$. Hence val(D)$=v^{\ast\ast}(b)\leq v(b)=$val(P).

\begin{proposition}
\label{cor-goz} For the problems \emph{(P)} and \emph{(D)} above, one has:
$\partial v(b)\neq\emptyset\Leftrightarrow$ [\emph{val(P)}$={}$\emph{val(D)}%
$\in\mathbb{R}$ and \emph{(D)} has optimal solutions]; moreover, if $\partial
v(b)\neq\emptyset$, then $\partial v(b)$ is the set \emph{Sol(D)} of the
optimal solutions of the dual problem \emph{(D)}.
\end{proposition}

\emph{Proof.} Assume that $\partial v(b)\neq\emptyset$ and take $y^{\ast}%
\in\partial v(b)$. Then $v(b)\in\mathbb{R}$ and $v(b)+v^{\ast}(y^{\ast
})=\left\langle b,y^{\ast}\right\rangle $ by (\ref{r-goz2}). Consequently,
val(P)$=v(b)=\left\langle b,y^{\ast}\right\rangle -v^{\ast}(y^{\ast})\leq{}%
$val(D), and so val(P)$={}$val(D)$;$ hence $y^{\ast}\in$Sol(D). Conversely,
assume that $[v(b)=]$ val(P)$={}$val(D)$\in\mathbb{R}$ and (D) has an optimal
solution $y^{\ast}$. Then $(\mathbb{R}\ni)$ $v(b)=\left\langle b,y^{\ast
}\right\rangle -v^{\ast}(y^{\ast})$, and so $y^{\ast}\in\partial v(b);$ hence
Sol(D)$\subset\partial v(b)$. Therefore, Sol(D)$=\partial v(b)$ whenever
$\partial v(b)\neq\emptyset$. \hfill$\square$

\medskip

Proposition \ref{cor-goz} is essentially \cite[Prop.\ 2.5]{Sha01}; its first
part is established in \cite[Th.\ 1]{GrOsZa02}. In this context it is worth
recalling \cite[Rem.\ 1, p.\ 274]{GrOsZa02}:

\medskip

Q1\label{Q1} -- \textquotedblleft REMARK 1. The Lipschitz property of the
value function $v$ in the assignment model ensures that the set of dual
solutions coincides with the subdifferential of the value function.
\emph{This}, of course, \emph{need not be true in more general economic
models}." (Our emphasis.)

\smallskip

The second part of Proposition \ref{cor-goz} shows that, for the general conic
linear programming problem (P), the set of solutions of the dual problem (D)
coincides with the subdifferential of the value function whenever the latter
is nonempty, not only for \textquotedblleft the value function $v$ in the
assignment model".

\medskip

In \cite[p.\ 266]{GrOsZa02}, one says \textquotedblleft we will consider only
LP problems for which the value function is proper"; hence, in \cite{GrOsZa02}%
, $v$ is a proper sublinear function with $\operatorname*{dom}v=A(P)-Q$. Of
course, in this case (that is, $v$ is proper), if $\dim Y<\infty$ then $v$ is
sub\-differentiable on $\operatorname*{icr}(\operatorname*{dom}v)$
$(\neq\emptyset)$, and so $\emptyset\neq\partial v(y)$ [$\subset\partial v(0)$
by (\ref{r-goz18c})] for $y\in\operatorname*{icr}(\operatorname*{dom}v)$,
whence $c^{\ast}\in P^{+}+A^{\ast}(Q^{+})$.

\medskip The following statement seems to be an important result from
\cite{GrOsZa02}, even if it is not mentioned explicitly throughout this article.

\medskip

Q2\label{Q2} -- \textquotedblleft LEMMA 1. \emph{If the value function} $v$
for a linear programming problem in standard form on ordered normed linear
spaces \emph{is proper}, \emph{then} $v$ \emph{is a lower semicontinuous}
extended real-valued (convex and homogeneous) \emph{function.}" (Our emphasis.)

\medskip

We shall see in the next section that \cite[Lem.\ 1]{GrOsZa02} is false even
in finite dimensional spaces.

\medskip Immediately after the proof of \cite[Lem.\ 1]{GrOsZa02} one finds the
following text:

\medskip

Q3\label{Q3} -- \textquotedblleft As promised in the Introduction, we can
quickly derive the Duffin--Karlovitz [3]\footnote{This is our reference
\cite{DufKar65}.} and the Charnes--Cooper--Kortanek [1, 2] no-gap theorems
from Theorem 1. \emph{Let }$X$\emph{ and }$Y$\emph{ be ordered normed linear
spaces and consider a linear programming problem in standard form with data
}$A$, $b$, $c^{\ast}$. The Duffin--Karlovitz theorem asserts that \emph{if the
positive cone} $Y_{+}$ \emph{has non-empty interior, if there is a feasible
solution} $\widehat{x}$ \emph{for the primal problem such that} $\widehat
{x}\geq0$ \emph{and} $A\widehat{x}-b$ \emph{is in the interior of} $Y_{+}$,
\emph{and if the value of the primal is finite, then the dual problem has a
solution and there is no gap.}\footnote{We did not find an assertion
equivalent to the mentioned \textquotedblleft Duffin--Karlovitz theorem" in
\cite{DufKar65}; in fact, in \cite[p.\ 123]{DufKar65}, one says:
\textquotedblleft This theory makes very little use of topology so it is more
like the theory of finite linear programming than like the theories given in
[2] and [5]. The desirability of omitting topological considerations is
emphasized by the paper of Charnes, Cooper and Kortanek. (However, in another
paper [4] a topological approach to this and similar problems will be
treated.)".} Since $\widehat{x}$ is feasible for $b$, $b$ is in the domain of
the value function $v$ and, hence, $v(b)<+\infty$. By hypothesis there is an
open ball $U$ around $A\widehat{x}-b$ within $F_{+}$. Hence, $\widehat{x}$ is
feasible for $b+u$ for all $u\in U$; viz.\ $b$ is an interior point of the
domain of $v$. The hypothesis that the value of the primal is finite is just
the statement that $v(b)>-\infty$. Consequently, $v$ is subdifferentiable at
$b$ and we conclude by Theorem 1 that there exists a dual solution and that
there is no gap, as asserted." (Our emphasis.)

\smallskip

We have to understand that the statement of the Duffin--Karlovitz theorem is
given by the emphasized text from Q3. Let us analyze the given proof. We agree
with the facts that $b+U\subset\operatorname*{dom}v$ [and so $b\in
\operatorname*{int}(\operatorname*{dom}v)$] and $v(b)\in\mathbb{R}$. From
this, without providing a motivation, one concludes that $\partial
v(b)\neq\emptyset$. Even if not mentioned, it is true that $v(y)\in\mathbb{R}$
for every $y\in\operatorname*{dom}v$, and so $v$ is proper. Using
Lemma\ 1,\footnote{The string \textquotedblleft Lemma 1" appears only once in
\cite{GrOsZa02}, more precisely in the text from Q2.} it follows that $v$ is
lower semi\-continuous. Why is $\partial v(b)$ nonempty? Without being
mentioned explicitly, probably, one uses the following statement from
\cite[p.\ 266]{GrOsZa02}:

\medskip

Q4\label{Q4} -- \textquotedblleft Suppose that $f:Y\rightarrow\mathbb{R}%
\cup\{+\infty\}$ is a proper convex function defined on the normed linear
space $Y$ and that $b\in\operatorname*{dom}f$. Each of the following
conditions implies the next and the last is equivalent to the
subdifferentiability of $f$ at $b$.

1. $f$ is lower semicontinuous and $b$ is an interior point of
$\operatorname*{dom}f$;

2. $f$ is locally Lipschitz at the point $b$, i.e.\ there exists $\delta>0$
such that $f$ is Lipschitz on $\operatorname*{dom}f\cap B(b;\delta)$;

3. $f$ has bounded steepness at the point $b$, i.e.\ the quotients
$(f(b)-f(y))/\left\Vert y-b\right\Vert $ are bounded above."

\smallskip

So, using \cite[Lem.\ 1]{GrOsZa02} and the above implications 1.\ $\Rightarrow
$ 2.\ $\Rightarrow$ 3., one gets $\partial v(b)\neq\emptyset$. Is the
implication 1.\ $\Rightarrow$ 2. true?

\smallskip

On page 267 of \cite{GrOsZa02} one mentions:

\medskip Q5\label{Q5} -- \textquotedblleft It is well-known (see Phelps
[9]\footnote{This is our reference \cite{Phe93}.}) that an extended
real-valued proper lower semicontinuous convex function is locally Lipschitz
and locally bounded on the interior of its domain."

\smallskip We did not succeed to find this assertion in (our reference)
\cite{Phe93}, but we found the following two related results:

\medskip

Q6\label{Q6} -- \textquotedblleft Proposition 1.6. If the convex function $f$
is continuous at $x_{0}\in D$, then it is locally Lipschitzian at $x_{0}$,
that is, there exist $M>0$ and $\delta>0$ such that $B(x_{0};\delta)\subset D$
and $\left\vert f(x)-f(y)\right\vert \leq M\left\Vert x-y\right\Vert $
whenever $x,y\in B(x_{0};\delta)$.

\smallskip

Proposition 3.3. Suppose that $f$ is a proper lower semicontinuous convex
function on a Banach space $E$ and that $D=\operatorname*{int}%
\operatorname*{dom}(f)$ is nonempty; then $f$ is continuous on $D$."

\smallskip

In fact, this version of \textquotedblleft the Duffin--Karlovitz theorem" is
contained in \cite[Th.\ 3]{Kre61} because $P:=X_{+}$ and $Q:=Y_{+}$ are
tacitly assumed to be closed in \cite{GrOsZa02}. Note that \textquotedblleft
the Duffin--Karlovitz theorem" is true even for $P$ and $Q$ not necessarily
closed; for this one could apply \cite[Th.\ 2.7.1]{Zal02} under its condition
(iii) for $\Phi:X\times Y\rightarrow\overline{\mathbb{R}}$ defined by
$\Phi(x,y):=\left\langle x,c^{\ast}\right\rangle +\iota_{P}(x)+\iota
_{Q}(Ax-y-b)$.

\section{Two examples}

In the sequel, the topological duals of Hilbert spaces (including
$\mathbb{R}^{m}$ with $m\in\mathbb{N}^{\ast}$) are identified with themselves
using Riesz' theorem.

The next example shows that \cite[Lem.\ 1]{GrOsZa02} is not true even in
finite dimensional spaces.

\begin{example}
\label{ex1} Consider $X:=\mathbb{R}^{2}$, $Y:=\mathbb{R}^{2}\times\mathbb{R}$,
$A:X\rightarrow Y$ with $A(x_{1},x_{2}):=(x_{1},x_{2},0)$, $c^{\ast}%
:=(0,1)\in\mathbb{R}^{2}$,
\[
P:=\mathbb{R}\times\mathbb{R}_{+},\quad Q:=\left\{  (y_{1},y_{2},y_{3})\in
Y\mid y_{1},y_{3}\in\mathbb{R}_{+},\ (y_{2})^{2}\leq2y_{1}y_{3}\right\}  ,
\]
the conic linear programming problem

\smallskip\emph{(P)} \ minimize\ $\left\langle x,c^{\ast}\right\rangle $
\ s.t. $x\in P$, $Ax-b\in Q$,

\smallskip\noindent and $v:Y\rightarrow\overline{\mathbb{R}}$ defined by
$v(y):=\inf\{\left\langle x,c^{\ast}\right\rangle \mid x\in P$, $Ax-b\in Q\}$.

Then $v(y)=y_{2}$ if $y:=(y_{1},y_{2},y_{3})\in\mathbb{R}\times\mathbb{R}%
_{+}\times\{0\}$, $v(y)=0$ for $y\in\mathbb{R}\times\mathbb{R}\times
(-\mathbb{P})$, and $v(y)=\infty$ elsewhere. Consequently, $v$ is not lower
semi\-continuous at any $y\in\mathbb{R}\times\mathbb{P}\times\{0\}$.
\end{example}

\emph{Proof.} Observe that $\left\langle x,c^{\ast}\right\rangle =x_{2}\geq0$
for $x\in P$, and so $v(y)\geq0$ for every $y\in\operatorname*{dom}v=A(P)-Q$.
Take $y\in Y$ and $x\in P$; if $Ax-y=(x_{1}-y_{1},x_{2}-y_{2},-y_{3})\in Q$,
then $x_{2}\geq0$, $x_{1}\geq y_{1}$, $y_{3}\leq0$ and $(x_{2}-y_{2})^{2}%
\leq-2y_{3}(x_{1}-y_{1})$.

Hence $y\notin\operatorname*{dom}v$ for $y_{3}>0$, and so $v(y)=\infty$. Take
$y_{3}=0$; then necessarily $x_{2}=y_{2}$. Hence $y\notin\operatorname*{dom}v$
for $y_{2}<0$, and so $v(y)=\infty$. If $y_{2}\geq0$ then $x_{2}=y_{2};$
taking $x_{1}=y_{1}$, $x=(x_{1},x_{2})$ is feasible for (P), and so one
obtains $v(y)=y_{2}$ in this case.

Take now $y_{3}<0$; then $x=(x_{1},0)$ with $x_{1}=y_{1}-\tfrac{1}{2}%
(y_{2})^{2}/y_{3}$ is feasible for (P), and so $v(y)=0$ in this case.

Consequently $v(y)=y_{2}$ if $y\in\mathbb{R}\times\mathbb{R}_{+}\times\{0\}$,
$v(y)=0$ for $y\in\mathbb{R}\times\mathbb{R}\times(-\mathbb{P})$, and
$v(y)=\infty$ elsewhere; hence,
\[
\operatorname*{dom}v=\left(  \mathbb{R}\times\mathbb{R}_{+}\times\{0\}\right)
\cup\left(  \mathbb{R}\times\mathbb{R}\times(-\mathbb{P})\right)  .
\]

Clearly, $v$ is convex (in fact sublinear), but $v$ is not l.s.c.\ at any
$y\in\mathbb{R}^{3}$ with $y_{2}>0$ and $y_{3}=0$; indeed, in this case,
$\operatorname*{dom}v\ni\zeta_{n}:=(y_{1},y_{2},-1/n)\rightarrow
y:=(y_{1},y_{2},0)$ and $v(\zeta_{n})=0\rightarrow0<y_{2}=v(y)$.
\hfill$\square$

\medskip

The next example is an adaptation of \cite[Examp.\ 2.3]{Zal08a} to the present
context. It shows that the value function $v$ can be proper and not lower
semi\-continuous at each $y\in\operatorname*{dom}v$.

\begin{example}
\label{ex5}Let $X$ be a separable infinite-dimensional real Hilbert space with
the orthonormal basis $(e_{n})_{n\geq1}$ and
\[
P:=\big\{
\sum\nolimits_{n\geq1}\lambda_{n}z_{n}\mid(\lambda_{n})\in(\ell_{2}%
)_{+}\big\}  \subset X,\quad c^{\ast}:=\sum\nolimits_{n\geq1}\eta_{n}e_{2n},
\]
with $z_{n}:=\eta_{n}e_{2n-1}-\mu_{n}e_{2n}$, where $\eta_{n},\mu_{n}\in
{}]0,1[$ are such that $\eta_{n}^{2}+\mu_{n}^{2}=1$ for every $n\geq1$ and
$(\eta_{n})_{n\geq1}\in\ell_{2}$. Consider $\Pr_{L}:X\rightarrow X$ the
orthogonal projection on
\[
L:=\overline{\operatorname*{span}}\left\{  e_{2n-1}\mid n\geq1\right\}
=\big\{\sum\nolimits_{n\geq1}\lambda_{n}e_{2n-1}\mid(\lambda_{n})\in\ell
_{2}\big\},
\]
$A:X\rightarrow Y:=X$ defined by $Ax:=\Pr_{L}x$, the conic linear programming problems

\smallskip\emph{(P)} \ minimize\ $\left\langle x,c^{\ast}\right\rangle $
\ s.t. $x\in P$, $Ax=b$,

\smallskip\noindent and $v:Y\rightarrow\overline{\mathbb{R}}$ defined by
$v(y):=\inf\{\left\langle x,c^{\ast}\right\rangle \mid x\in P$, $Ax=y\}$ for
$y\in Y$.

Then
\begin{gather}
\operatorname*{dom}v=\big\{y:=\sum\nolimits_{n\geq1}\gamma_{n}e_{2n-1}%
\mid\big(\left(  \eta_{n}\right)  ^{-1}\gamma_{n}\big)_{n\geq1}\in(\ell
_{2})_{+}\big\}\subset L,\label{r-goz21a}\\
v(y)=-\sum\nolimits_{n\geq1}\mu_{n}\gamma_{n}\leq0=v(0),\quad\forall
y:=\sum\nolimits_{n\geq1}\gamma_{n}e_{2n-1}\in\operatorname*{dom}v.
\label{r-goz21b}%
\end{gather}
More precisely, \emph{(P)} has a unique feasible solution (hence a unique
optimal solution) for every $b\in\operatorname*{dom}v$, and so $v$ is a proper
sublinear function. Moreover, the dual problem \emph{(D)} has not feasible
solutions for every $b\in Y$, proving so that $v$ is not l.s.c.\ at any
$b\in\operatorname*{dom}v$; in particular, $\partial v(b)=\emptyset$ for every
$b\in Y$.
\end{example}

Proof. Note that $\left\langle z_{n},z_{m}\right\rangle =\delta_{nm}$ for
$n,m\geq1$ ($\delta_{nm}$ being the Kronecker's symbols). Clearly, $P$ is a
closed convex cone. Consider $Q:=\{0\}\subset Y$; then $Q^{+}=Y$. Because
$\Pr_{L}=(\Pr_{L})^{\ast}$, one obtains that
\[
\operatorname*{dom}v=A(P)-Q=\Pr\nolimits_{L}(P),\quad\operatorname*{dom}%
v^{\ast}=Q^{+}\cap(A^{\ast})^{-1}(c^{\ast}-P^{+})=\Pr\nolimits_{L}%
^{-1}(c^{\ast}-P^{+}).
\]

Consider now $y\in\operatorname*{dom}v$ $(\subset L)$; hence $y=\sum_{n\geq
1}\gamma_{n}e_{2n-1}$ with $(\gamma_{n})\in\ell_{2}$, and there exists $x\in
P$ such that $y=Ax$; hence $x=\sum\nolimits_{n\geq1}\lambda_{n}z_{n}%
=\sum\nolimits_{n\geq1}\lambda_{n}(\eta_{n}e_{2n-1}-\mu_{n}e_{2n})$ for some
$(\lambda_{n})\in(\ell_{2})_{+}$. Therefore, $\gamma_{n}=\lambda_{n}\eta
_{n}\geq0$, whence $\lambda_{n}=\gamma_{n}/\eta_{n}$, for $n\geq1$; this shows
that the set $\{x\in P\mid y=Ax\}$ is a singleton $\{x_{y}\}$ for (every)
$y\in\operatorname*{dom}v$ and so
\[
v(y)=\left\langle c^{\ast},x_{y}\right\rangle =\sum\nolimits_{n\geq1}\eta
_{n}(-\lambda_{n}\mu_{n})=\sum\nolimits_{n\geq1}\eta_{n}(-\mu_{n}\cdot
\gamma_{n}/\eta_{n})=-\sum\nolimits_{n\geq1}\mu_{n}\gamma_{n}.
\]

Consequently, (\ref{r-goz21a}) and (\ref{r-goz21b}) hold.

Assume that $x\in\operatorname*{dom}v^{\ast}$; then there exists $(\lambda
_{n})\in\ell_{2}$ such that $u:=\Pr_{L}x=\sum\nolimits_{n\geq1}\lambda
_{n}e_{2n-1}\in c^{\ast}-P^{+}$, that is, $c^{\ast}-u\in P^{+}$; hence
\[
0\leq\left\langle c^{\ast}-u,z_{k}\right\rangle =-\lambda_{k}\eta_{k}-\eta
_{k}\mu_{k}=-\eta_{k}(\lambda_{k}+\mu_{k})\quad\forall k\geq1.
\]
It follows that $\lambda_{k}+\mu_{k}\leq0$ for all $k\geq1$, contradicting the
fact that $\lambda_{n}\rightarrow0$ and $\mu_{n}\rightarrow1$. Hence
$\operatorname*{dom}v^{\ast}=\emptyset$, and so (D$_{y}$) has not feasible
solutions for every $y\in Y$. Consequently, $v^{\ast\ast}=-\infty$, proving so
that $\overline{v}(y)=-\infty$ for every $y\in\operatorname*{dom}v$, and so
$v$ is not l.s.c.\ at any $y\in\operatorname*{dom}v$. \hfill$\square$

\section{On Kretschmer's gap example in linear programming}

In \cite[pp.\ 230, 231]{Kre61}, Kretschmer considers $Y:=L^{2}:=L^{2}[0,1]$
(with respect to the Lebesgue measure $\mu$ on $[0,1]$) endowed with the usual
inner product and ordered by $Q:=L_{+}^{2}$, as well as $X:=L^{2}%
\times\mathbb{R}$ endowed with the inner product defined by $\left\langle
(x,r),(x^{\prime},r^{\prime})\right\rangle :=\left\langle x,x^{\prime
}\right\rangle +rr^{\prime}$ and ordered by $P:=L_{+}^{2}\times\mathbb{R}_{+}%
$; obviously, $P^{+}=P$ and $Q^{+}=Q$. Moreover, one takes $A:X\rightarrow Y$
with $A(x,r):=y+re_{0}$, with $y(t):=\int_{t}^{1}x(s)ds$ for $t\in\lbrack0,1]$
and $e_{0}\in L^{2}$ with $e_{0}(t):=1$ for $t\in\lbrack0,1]$. Furthermore,
$A$ is a continuous linear operator and $A^{\ast}$ $(:Y\rightarrow X)$ is
given by $A^{\ast}y=(x,r)\in X$ with $x(t):=\int_{0}^{t}y(s)ds$ for
$t\in\lbrack0,1]$ and $r=\int_{0}^{1}y(s)ds$.

\smallskip Let $c^{\ast}:=c_{\alpha}^{\ast}:X\rightarrow\mathbb{R}$ be defined
by $c^{\ast}(x,r):=\int_{0}^{1}tx(t)dt+\alpha r=\left\langle (x,r),(e_{1}%
,\alpha)\right\rangle $, where $\alpha\in\mathbb{R}_{+}$ and $e_{1}(t):=t$ for
$t\in\lbrack0,1]$; clearly, $c^{\ast}\in P^{+}$.

\smallskip Consider the problem (P):=(P$_{\alpha}$) and its dual
(D):=(D$_{\alpha}$) defined by

\smallskip(P) \ minimize $\left\langle (x,r),c^{\ast}\right\rangle $
\ s.t.\ \ $(x,r)\geq0$, $A(x,r)-b\geq0$,

\smallskip(D) \ maximize $\left\langle z,b\right\rangle $ \ s.t.\ \
$z\geq0$, $A^{\ast}z\leq c^{\ast}$,

\smallskip\noindent as well as the value function%
\[
v:=v_{\alpha}:Y\rightarrow\overline{\mathbb{R}},\quad v_{\alpha}%
(y):=\inf\{\left\langle (x,r),c_{\alpha}^{\ast}\right\rangle \mid(x,r)\in
F(y)\},
\]
where
\[
F(y):=\{(x,r)\in X\mid(x,r)\geq0,\ A(x,r)\geq y\}
\]
is the feasible set of the problem (P); notice that $F(y)$ is the same for all
$\alpha\in\mathbb{P}$. Clearly, $v(y)\geq0=v(0)$ for $y\in Y$ because
$c^{\ast}\in P^{+}$; hence $0\in\partial v(0)$. In fact, by (\ref{r-goz1}),
(\ref{r-goz18a}) and (\ref{r-goz18b}), one has
\[
\partial v(0)=Q^{+}\cap(A^{\ast})^{-1}(c^{\ast}-P^{+})=Q\cap(A^{\ast}%
)^{-1}(c^{\ast}-P)\text{ \ and\ \ }v^{\ast\ast}=\overline{v}.
\]

Let us denote by $\mathcal{A}$ the class of measurable subsets of $[0,1]$.
Without loss of generality we assume that $y(t)\in\mathbb{R}$ for all $y\in
L^{2}$ and $t\in\lbrack0,1]$. For $y\in L^{2}$ and $\gamma\in\mathbb{R}$ we
set $[y\geq\gamma]:=\{t\in\lbrack0,1]\mid y(t)\geq\gamma\}$ and $y_{+}%
:=\max\{y,0\}$; clearly, $[y\geq\gamma]\in\mathcal{A}$ and $y_{+}\in L_{+}%
^{2}$. Moreover, the \emph{characteristic function} of $E\subset\lbrack0,1]$
is the function $\chi_{E}:[0,1]\rightarrow\mathbb{R}$ defined by $\chi
_{E}(t):=1$ for $t\in E$ and $\chi_{E}(t):=0$ for $t\in\lbrack0,1]\setminus E$.

\begin{lemma}
\label{l-Z1}The following assertions hold:

\emph{(i)} One has
\[
\operatorname*{dom}v=\{y\in Y\mid\operatorname*{ess}\sup y<\infty\}=\{y\in
L^{2}\mid y_{+}\in L^{\infty}\}; \label{r-goz4}%
\]
in particular $L^{\infty}\subset\operatorname*{dom}v$, and so
$\operatorname*{cl}(\operatorname*{dom}v)=Y$.

\emph{(ii)} Let $A\in\mathcal{A}$ be such that $\beta:=\mu(A)>0$. Then there
exists a sequence $(A_{n})_{n\geq1}\subset\mathcal{A}$ such that
$A=\cup_{n\geq1}A_{n}$, $A_{n}\cap A_{m}=\emptyset$ for $n\neq m$ and
$\mu(A_{n})=2^{-n}\beta$ for $n\geq1$.

\emph{(iii)} Let $A$ and $(A_{n})_{n\geq1}$ be as in \emph{(ii)} and consider
\[
\widetilde{y}_{n}:=\sum\nolimits_{k=1}^{n}2^{k/4}\chi_{A_{k}}\geq
0\ \ (n\geq1),\quad\widetilde{y}:=\sup\nolimits_{n\geq1}\widetilde{y}_{n}%
\geq0. \label{r-goz12}%
\]
Then $\operatorname*{ess}\sup\widetilde{y}_{n}=2^{n/4}\rightarrow\infty$,
$\widetilde{y}\in L^{2}$, $\left\Vert \widetilde{y}_{n}\right\Vert <\left\Vert
\widetilde{y}\right\Vert =\left[  \beta(\sqrt{2}+1)\right]  ^{1/2}$ and
$\left\Vert \widetilde{y}_{n}-\widetilde{y}\right\Vert \rightarrow0$;
consequently $\widetilde{y}_{n}\in L_{+}^{\infty}\subset L_{+}^{2}$ for
$n\geq1$ and $\widetilde{y}\in L_{+}^{2}\setminus L^{\infty}$.
\end{lemma}

Proof. (i) Let $y\in\operatorname*{dom}v$; then there exists $(x,r)\in P$ such
that $y\leq A(x,r)$, and so $y(t)\leq\int_{t}^{1}x(s)ds+r\leq\int_{0}%
^{1}x(s)ds+r\leq\left\Vert x\right\Vert +r$ for$\ t\in\lbrack0,1]$, whence
$\operatorname*{ess}\sup y<\infty$. Conversely, if $y\in L^{2}$ is such that
$r:=\operatorname*{ess}\sup y<\infty$, then $y\leq_{Q}re_{0}\leq_{Q}r_{+}%
e_{0}=A(0,r_{+})$, where $r_{+}:=\max\{0,r\}$, and so $v(y)\leq\left\langle
(0,r_{+}),c^{\ast}\right\rangle =\alpha r_{+}$; hence $y\in\operatorname*{dom}%
v$.

(ii) Because $\mu$ has not atoms, there exists $A_{1}\subset A$ such that
$A_{1}\in\mathcal{A}$ and $\mu(A_{1})=2^{-1}\beta$ $(\in{}]0,\mu(A)[)$. Hence
$A_{1}^{\prime}:=A\setminus A_{1}\in\mathcal{A}$ and $\mu(A_{1}^{\prime}%
)=\mu(A)-\mu(A_{1})>2^{-2}\beta$, and so there exists $A_{2}\subset
A_{1}^{\prime}$ such that $A_{2}\in\mathcal{A}$ and $\mu(A_{2})=2^{-2}\beta$
$(\in{}]0,\mu(A_{1}^{\prime})[)$; clearly, $A_{1}\cap A_{2}=\emptyset$.
Continuing in the same way we get the sequence $(A_{n})_{n\geq1}%
\subset\mathcal{A}$ with $A_{n}\subset A$, $\mu(A_{n})=2^{-n}\beta$ and
$A_{n}\cap A_{m}=\emptyset$ for $n,m\in\mathbb{N}^{\ast}$ with $n\neq m$.
Setting $A^{\prime}:=\cup_{n\geq1}A_{n}$ one has that $A^{\prime}\subset A$
and $\mu(A^{\prime})=\sum_{n\geq1}\mu(A_{n})=\sum_{n\geq1}2^{-n}\beta
=\beta=\mu(A)$, and so $\mu(A\setminus A^{\prime})=0$. Replacing $A_{1}$ with
$A_{1}\cup(A\setminus A^{\prime})$, it follows that the sequence
$(A_{n})_{n\geq1}$ has the desired properties.

(iii) Because $A_{k}\in\mathcal{A}$ for $k\geq1$, one has that $\widetilde
{y}_{n}$ is measurable for $n\geq1$. Moreover, because $A_{n}\cap
A_{m}=\emptyset$ for $n\neq m$ and $\mu(A_{n})=2^{-n}\beta>0$ for $n\geq1$, it
is clear that $\widetilde{y}_{n}(t)=2^{n/4}\geq\widetilde{y}_{n}(s)$ for all
$t\in A_{n}$ and $s\in\lbrack0,1]$, and so $\operatorname*{ess}\sup
\widetilde{y}_{n}=\sup\widetilde{y}_{n}=2^{n/4}$ for $n\geq1$. Hence
$(\widetilde{y}_{n})_{n\geq1}\subset L_{+}^{\infty}\subset L_{+}^{2}$.

Clearly, $\widetilde{y}$ is measurable and $\widetilde{y}(t)=\lim
_{n\rightarrow\infty}\widetilde{y}_{n}(t)=\sum\nolimits_{k=1}^{\infty}%
2^{k/4}\chi_{A_{k}}(t)$ for $t\in\lbrack0,1]$. Because $\widetilde{y}_{n}%
\leq\widetilde{y}$, one has $2^{n/4}=\operatorname*{ess}\sup\widetilde{y}%
_{n}\leq\operatorname*{ess}\sup\widetilde{y}$ for $n\geq1$, and so
$\operatorname*{ess}\sup\widetilde{y}=\infty$; hence $\widetilde{y}\notin
L^{\infty}$. On the other hand, for $t\in\lbrack0,1]$ one has%
\[
(\widetilde{y}_{n})^{2}(t)=\sum\nolimits_{k=1}^{n}2^{k/2}\chi_{A_{k}%
}(t)\ \ (n\geq1),\quad\widetilde{y}^{2}(t)=\sum\nolimits_{k=1}^{\infty}%
2^{k/2}\chi_{A_{k}}(t),
\]
and so
\begin{gather*}
\left\Vert \widetilde{y}_{n}\right\Vert ^{2}=\sum\nolimits_{k=1}^{n}2^{k/2}%
\mu(A_{k})=\beta\sum\nolimits_{k=1}^{n}2^{-k/2}=(\sqrt{2}+1)(1-2^{-n/2}%
)\beta,\\
\left\Vert \widetilde{y}\right\Vert ^{2}=\beta\sum\nolimits_{k=1}^{\infty
}2^{-k/2}=(\sqrt{2}+1)\beta,\quad\left\Vert \widetilde{y}_{n}-\widetilde
{y}\right\Vert ^{2}=2^{-n/2}(\sqrt{2}+1)\beta.
\end{gather*}
Therefore, $\widetilde{y}\in L_{+}^{2}\setminus L^{\infty}$, $\left\Vert
\widetilde{y}_{n}\right\Vert <\left\Vert \widetilde{y}\right\Vert =\left[
\beta(\sqrt{2}+1)\right]  ^{1/2}$ and $\left\Vert \widetilde{y}_{n}%
-\widetilde{y}\right\Vert \rightarrow0$. \hfill$\square$

\begin{proposition}
\label{p-KZ}Assume that $\alpha>0$. Then for every $y\in\operatorname*{dom}v$
and every $\rho>0$, $v$ is not bounded on $B(y,\rho)\cap\operatorname*{dom}v$;
in particular, $v|_{\operatorname*{dom}v}$ is not continuous at each
$y\in\operatorname*{dom}v$.
\end{proposition}

Proof. Consider $y\in\operatorname*{dom}v$ and $\rho>0$, as well as
$0<\eta_{0}<\eta_{1}<\min\{1,\alpha\}$ and $B:=[0,\eta_{0}]\cup\lbrack\eta
_{1},1]$; hence $\operatorname*{ess}\sup y<\infty$. For each $k\in
\mathbb{N}^{\ast}$, set $E_{k}:=[y\geq-k]$; because $[0,1]=\cup_{k\geq1}E_{k}$
and $E_{k}\subset E_{k+1}$ for $k\geq1$, it follows that $\mu(E_{k}%
)\rightarrow\mu([0,1])=1$. Consider $k_{0}\geq1$ such that $\mu(E_{k_{0}%
})>\eta_{0}+1-\eta_{1}=\mu(B)$, and set $\gamma:=-k_{0}$, $A:=E_{k_{0}%
}\setminus B=E_{k_{0}}\cap{}]\eta_{0},\eta_{1}[{}\subset{}]\eta_{0},\eta_{1}[$
and $\beta:=\mu(A)$. Clearly, $E_{k_{0}}\subset A\cup B$, and so $\mu
(B)<\mu(E_{k_{0}})\leq\mu(A)+\mu(B)=\beta+\mu(B)$, whence $\beta>0$; set also
$\delta:=\left[  \beta(\sqrt{2}+1)\right]  ^{1/2}>0$.

Consider now the sets $A_{n}$ and the functions $\widetilde{y}_{n}$ for
$n\geq1$ provided by assertions (ii) and (iii) of Lemma \ref{l-Z1}, as well as
$\widetilde{y}:=\sup_{n\geq1}\widetilde{y}_{n}$; hence
\[
L_{+}^{\infty}\ni\widetilde{y}_{n}\rightarrow^{\left\Vert \cdot\right\Vert
}\widetilde{y}\in L_{+}^{2}\setminus L^{\infty}\ \ \text{and\ \ }\left[
\forall n\geq1:\left\Vert \widetilde{y}_{n}\right\Vert <\left\Vert
\widetilde{y}\right\Vert =\delta\right]  .
\]
Consider also $0<\varepsilon<\rho/\delta$ and set $y_{n}:=y+\varepsilon
\widetilde{y}_{n}$ for $n\in\mathbb{N}^{\ast}$; clearly $y_{n}\in L^{2}$ and
$\operatorname*{ess}\sup y_{n}\leq\operatorname*{ess}\sup y+\varepsilon
\operatorname*{ess}\sup\widetilde{y}_{n}<\infty$, whence $y_{n}\in
\operatorname*{dom}v$ and $\left\Vert y_{n}-y\right\Vert =\left\Vert
\varepsilon\widetilde{y}_{n}\right\Vert <\varepsilon\delta<\rho$ for $n\geq1$.
Moreover, $y_{n}\rightarrow^{\left\Vert \cdot\right\Vert }y+\varepsilon
\widetilde{y}\notin\operatorname*{dom}v$. Therefore,
\begin{equation}
(y_{n})_{n\geq1}\subset B(y,\rho)\cap\operatorname*{dom}v\ \ \text{and
\ }B(y,\rho)\cap(Y\setminus\operatorname*{dom}v)\neq\emptyset. \label{r-goz16}%
\end{equation}

Let $n\geq1$ be fixed and consider $(x,r)\in F(y_{n})$; hence%
\[
\textstyle
x\geq0,\ \ r\geq0,\ \ \text{and \ }\int_{t}^{1}x(s)ds+r\geq y(t)+\varepsilon
\widetilde{y}_{n}(t)\ \ \text{a.e.}\ t\in\lbrack0,1].\label{r-goz14}%
\]
Because $A_{n}\subset A=E_{n_{0}}\cap{}]\eta_{0},\eta_{1}[{}\subset\lbrack
y\geq\gamma])$, one has
\[
\textstyle
\int_{t}^{1}x(s)ds+r\geq y(t)+\varepsilon\widetilde{y}_{n}(t)\geq
\gamma+2^{n/4}\varepsilon\ \ \text{for a.e.}\ t\in A_{n},
\]
and so, for a.e.\ $t\in A_{n}$, one has%
\begin{align*}
\left\langle (x,r),c^{\ast}\right\rangle  &  =\textstyle\int_{0}%
^{1}sx(s)ds+\alpha r\geq\int_{t}^{1}sx(s)ds+\alpha r\geq t\int_{t}%
^{1}x(s)ds+\alpha r\\
&  \geq\eta_{0}(\gamma+2^{n/4}\varepsilon-r)+\alpha r=\eta_{0}(\gamma
+2^{n/4}\varepsilon)+r(\alpha-\eta_{0})\\
&  \geq\eta_{0}(\gamma+2^{n/4}\varepsilon)+r(\eta_{1}-\eta_{0})\geq\eta
_{0}(\gamma+2^{n/4}\varepsilon);
\end{align*}
hence $\left\langle (x,r),c^{\ast}\right\rangle \geq\eta_{0}(\gamma
+2^{n/4}\varepsilon)$. Because $(x,r)\in F(y_{n})$ is arbitrary, it follows
that $v(y_{n})\geq\eta_{0}(\gamma+2^{n/4}\varepsilon)$. Therefore,
$v(y_{n})\geq\eta_{0}(\gamma+2^{n/4}\rho/\delta)$ for every $n\geq1$, and so
$v(y_{n})\rightarrow\infty$. Taking into account (\ref{r-goz16}), it follows
that $v$ is not bounded on $B(y,\rho)\cap\operatorname*{dom}v$; moreover,
$y\notin\operatorname*{int}(\operatorname*{dom}v)$ because $B(y,\rho
)\cap(Y\setminus\operatorname*{dom}v)\neq\emptyset$ for every $\rho>0$,
proving that $\operatorname*{int}(\operatorname*{dom}v)=\emptyset$.
\hfill$\square$

\medskip

Observe that the case $\alpha=0$ is very special. Indeed, as seen in the proof
of Lemma \ref{l-Z1}(i), for $y\in\operatorname*{dom}v$, $(0,r_{+})\in F(y)$,
where $r:=\operatorname*{ess}\sup y$, and so $0\leq v(y)\leq\left\langle
(0,r_{+}),(e_{1},0)\right\rangle =0$. Hence $v(y)=0$ and the value $v(y)$ is
attained. Therefore, $v=\iota_{\operatorname*{dom}v}$. On the other hand, for
$y\in Y$, $z$ is feasible for the dual problem (D$_{y}$) if and only if
$z\geq0$, $\int_{0}^{t}z(s)ds\leq e_{1}(t)$ a.e.\ $t\in\lbrack0,1]$ and
$\int_{0}^{1}z(s)ds\leq\alpha=0$, and so $z=0$ is the only feasible (hence
optimal) solution of (D$_{y}$). Hence $v^{\ast\ast}(y)=0=\overline{v}(y)$ for
every $y\in Y$. Because $v=\iota_{\operatorname*{dom}v}$, one has
$\overline{v}=\iota_{\operatorname*{cl}(\operatorname*{dom}v)}$, confirming so
that $\operatorname*{cl}(\operatorname*{dom}v)=Y$.

\medskip

Taking $\alpha:=2$ and $b:=e_{0}$, one obtains \cite[Examp.\ 5.1]{Kre61}; this
is also considered in \cite[Examp.\ 1]{GrOsZa02}, as well as the one in which
$b:=b_{0}:=\chi_{\lbrack0,1/2]}$. The next two results are slight extensions
of those related to the \textquotedblleft modification" of \cite[Examp.\ 5.1]%
{Kre61} used in \cite[p.\ 270]{GrOsZa02}, the proofs using similar arguments
to those in \cite{GrOsZa02}.

\begin{proposition}
\label{ex4}Consider $\alpha\in\mathbb{P}$ and $b:=\chi_{I\cup J}$ with
$I:=[0,\delta]$, $J:=[\gamma,1]$, where $0\leq\delta\leq\gamma<1$. Then
\emph{val(P)}$=\alpha$, \emph{val(D)}$=\min\{1,\alpha\}$, and \emph{(P)},
\emph{(D)} have optimal solutions; moreover, \emph{val(P)} = \emph{val(D)}
$\Leftrightarrow$ $\alpha\in{}]0,1]$ $\Leftrightarrow$ $\partial v(\chi_{I\cup
J})\neq\emptyset$.
\end{proposition}

Proof. Clearly, if $(x,r)\in P$ is feasible then $\int_{t}^{1}x(s)ds+r\geq1$
a.e.\ $t\in\lbrack\gamma,1[;$ because $\lim_{[\gamma,1[\,\ni t\rightarrow
1}\int_{t}^{1}x(s)ds=0$, one gets $r\geq1$. Because $(0,1)$ is feasible for
(P), one has that $0$ is optimal solution for (P) and val(P)$=\alpha$.

Observe that for $z\geq0$ with $\int_{0}^{t}z(s)ds\leq t$ for $t\in
\lbrack0,1]$ one has $\int_{0}^{1}z(s)ds\leq1$, and so, when $z$ is feasible
for (D) one has $\int_{0}^{1}\chi_{\lbrack\gamma,1]}z=\int_{\gamma}^{1}%
z\leq\int_{0}^{1}z\leq\min\{1,\alpha\}$. Hence $0\leq$val(D)$\leq
\min\{1,\alpha\}=:\mu$. Take $\eta\in\lbrack\gamma,1[$ and $z:=\mu
(1-\eta)^{-1}\chi_{\lbrack\eta,1]}$ $(\geq0)$; then $\int_{0}^{t}z(s)ds=0$ for
$t\in\lbrack0,\eta]$ and $\int_{0}^{t}z(s)ds=\mu(1-\eta)^{-1}\int_{\eta}%
^{t}1ds=\mu\frac{t-\eta}{1-\eta}\leq\mu t\leq t$ for $t\in\lbrack\eta,1]$ and
so $z$ is feasible for (D). Moreover, $\int_{0}^{1}z(t)dt=\mu$, and so $z$ is
an optimal solution for (D), whence val(D)$=\min\{1,\alpha\}$. Consequently,
both problems have optimal solutions, and $\partial v(\chi_{\lbrack\gamma
,1]})\neq\emptyset$ if and only if $\alpha\in{}]0,1]$. \hfill$\square$

\medskip

Taking $\alpha:=2$ and $\delta:=\gamma=0$ one (re)obtains (as already
mentioned) the example from \cite[Examp.\ 5.1]{Kre61}, as well as the one from
\cite[p.\ 270]{GrOsZa02} and the conclusions from there, that is, both
problems have optimal solutions, but there is a (positive) duality gap.

\medskip

Consequently, the previous example shows not only that
$v|_{\operatorname*{dom}v}$ is not locally Lipschitz, but also that
$v|_{\operatorname*{dom}v}$ is not l.s.c.\ on its domain; therefore,
\cite[Examp.\ 1]{GrOsZa02} provides a counterexample to \cite[Lem.\ 1]%
{GrOsZa02}.

\begin{proposition}
\label{ex3}Take $\alpha\in\mathbb{P}$ and $b:=\chi_{\lbrack0,\delta]}$ with
$\delta\in{}]0,1[$. Then \emph{val(P)}$=$\emph{val(D)}$=\min\{\delta,\alpha\}$
and \emph{(D)} has optimal solutions; consequently, $\partial v(\chi
_{\lbrack0,\delta]})\neq\emptyset$. Furthermore, \emph{(P)} has optimal
solutions iff $\alpha\leq\delta$.
\end{proposition}

Proof. First observe that for $(x,r)\in P$, the following assertions are
equivalent: $(x,r)$ is feasible for (P); $(x\cdot\chi_{\lbrack\delta,1]},r)$
is feasible for (P); $\int_{\delta}^{1}x(s)ds+r\geq1$; $r\geq\big(1-\int
_{\delta}^{1}x(s)ds\big)_{+}$. Set%
\[
\textstyle F_{1}:=\big\{x\in L_{+}^{2}\mid\int_{\delta}^{1}x(s)ds\geq
1\big\},\quad F_{2}:=\big\{x\in L_{+}^{2}\mid\int_{\delta}^{1}x(s)ds\leq
1\big\}.
\]
Notice that $F_{1}\cap F_{2}\neq\emptyset$ and $0\in F_{2}$; moreover, $(x,0)$
is feasible when $x\in F_{1}$ and $\big(x,1-\int_{\delta}^{1}x(s)ds\big)$ is
feasible when $x\in F_{2}$. It follows that val(P)${}=\min\{v_{1},v_{2}\}$,
where $v_{1}:=\inf_{x\in F_{1}}\int_{0}^{1}tx(t)dt$ and
\begin{align*}
v_{2}:=  &  \inf_{x\in F_{2}}\big(\textstyle\int_{0}^{1}tx(t)dt+\alpha
-\alpha\int_{\delta}^{1}x(t)dt\big)=\inf_{x\in F_{2}}\big(\textstyle\alpha
+\int_{\delta}^{1}(t-\alpha)x(t)dt\big)\\
=  &  \,\alpha-\sup_{x\in F_{2}}\textstyle\int_{\delta}^{1}(\alpha
-t)x(t)dt\leq\alpha.
\end{align*}

For $x\in L_{+}^{2}$ and $t\in\lbrack\delta,1]$ one has $(\alpha
-t)x(t)\leq(\alpha-\delta)x(t)$, and so $\int_{\delta}^{1}(\alpha
-t)x(t)dt\leq(\alpha-\delta)\int_{\delta}^{1}x(t)dt$, with equality iff
$x\cdot\chi_{\lbrack\delta,1]}=0$. Assume that $x\in F_{2}$; for
$\alpha>\delta$ one has $\int_{\delta}^{1}(\alpha-t)x(t)dt\leq\alpha-\delta$,
while for $\alpha\leq\delta$ one has $\int_{\delta}^{1}(\alpha-t)x(t)dt\leq0$.
Therefore, $v_{2}\geq\delta$ if $\alpha\geq\delta$ and $v_{2}=\alpha$ if
$\alpha\leq\delta$, $v_{2}$ being attained for $x=0$ in the latter case.

In what concerns $v_{1}$, one has
\[
v_{1}=\inf_{x\in F_{1}}\big(  \textstyle
\int_{0}^{\delta}tx(t)dt+\int_{\delta}^{1}tx(t)dt\big)  =\inf_{x\in F_{1}%
}\textstyle\int_{\delta}^{1}tx(t)dt\geq\delta\inf_{x\in F_{1}}\textstyle\int
_{\delta}^{1}x(t)dt\geq\delta.
\]
For $\varepsilon\in{}]0,1-\delta\lbrack$ and $x:=\varepsilon^{-1}\chi
_{\lbrack\delta,\delta+\varepsilon]}$, one has $x\in F_{1}$ and $\int_{0}%
^{1}tx(t)dt=\varepsilon^{-1}\int_{\delta}^{\delta+\varepsilon}tdt=\delta
+\varepsilon/2$, and so $v_{1}=\delta$. Consequently, val(P)$=\min
\{\alpha,\delta\}$; moreover, if $\alpha>\delta$ then (P) has not optimal
solutions, and $x=0$ is solution of (P) if $\alpha\leq\delta$.

\smallskip

If $z$ is feasible for (D), then $\int_{0}^{\delta}z(t)dt\leq\delta$ and
$\int_{0}^{\delta}z(t)dt\leq\int_{0}^{1}z(t)dt\leq\alpha$, and so
val(D)$\leq\min\{\delta,\alpha\}$. Clearly, $z:=\chi_{\lbrack0,\min
\{\alpha,\delta\}]}$ is an optimal solution of (D), and so val(D)=$\min
\{\delta,\alpha\}$.

\smallskip Therefore, val(P) = val(D) $=\min\{\delta,\alpha\}$ and (D) has
optimal solutions; consequently, $\partial v(\chi_{\lbrack0,\delta]}%
)\neq\emptyset$. Furthermore, (P) has optimal solutions iff $\alpha\leq\delta
$. \hfill$\square$

\begin{corollary}
\label{c-goz} Let $\alpha\in{}]1,\infty\lbrack$ and $\delta\in{}]0,1[$, and
consider the problems

\smallskip\emph{(P$_{y}$)\ \ minimize} \ $\int_{0}^{1}tx(t)dt+\alpha r$ \ s.t.
\ $x\geq0$, $r\geq0$, $\int_{t}^{1}x(s)ds+r\geq y(t)$ a.e.\ $t\in\lbrack0,1]$,

\smallskip\emph{(D$_{y}$) \ maximize} \ $\int_{0}^{1}y(t)z(t)dt$ s.t
\ $z\geq0$, $\int_{0}^{t}z(s)ds\leq t$ a.e.\ $t\in\lbrack0,1]$, $\int_{0}%
^{1}z(s)ds\leq\alpha$.

\smallskip Then $\partial v(\chi_{\lbrack0,\delta]})\neq\emptyset$ and
$v|_{\operatorname*{dom}v}$ is not continuous at $\chi_{\lbrack0,\delta]}$.
\end{corollary}

Proof. By Proposition \ref{ex3} one has that $v(\chi_{\lbrack0,\delta
]})=\delta$ and $\partial v(\chi_{\lbrack0,\delta]})\neq\emptyset$, while from
Proposition \ref{ex4} one has that $v(\chi_{\lbrack0,\delta]\cup\lbrack
\gamma,1]})=\alpha$ for every $\gamma\in{}]\delta,1[$. Because $\Vert
\chi_{\lbrack0,\delta]\cup\lbrack\gamma,1]}-\chi_{\lbrack0,\delta]}\Vert
_{2}=\Vert\chi_{\lbrack\gamma,1]}\Vert_{2} =(1-\gamma)^{1/2}\rightarrow0$ for
$\gamma\rightarrow1$, confirming that $v|_{\operatorname*{dom}v}$ is not
continuous at $\chi_{\lbrack0,\delta]}$. \hfill$\square$

\medskip

In the paragraph before \cite[Examp.\ 1, p.\ 269]{GrOsZa02}, one says:

\smallskip Q7\label{Q7} -- \textquotedblleft We give an example of a convex
function which is subdifferentiable but not locally Lipschitz by exhibiting a
linear programming problem for which the value function has this property. The
example takes place in the Banach lattice $L^{2}[0,1]$ (a space for which the
positive cone has empty interior) and \emph{for which} $\operatorname*{dom}%
v\supseteq L^{2}[0,1]_{+}$". (Our emphasis.)

\smallskip This text is completed by the following ones from \cite[p.\ 270]%
{GrOsZa02}:

\smallskip Q8\label{Q8} -- \textquotedblleft On the other hand, we will
establish that $v$ is not locally Lipschitz at $b_{0}$; in fact, $v$ is not
even continuous there (\emph{or anywhere})." (Our emphasis.)

\smallskip Q9\label{Q9} -- \textquotedblleft\emph{Similar perturbations} show
that $v$ is not continuous anywhere on $L^{2}[0,1]_{+}$. (\emph{Of course,
}$v$\emph{ is lower semicontinuous}.)" (Our emphasis.)

\medskip

As seen in Lemma \ref{l-Z1}, one has $L_{+}^{\infty}\subset L^{\infty}%
\subset\operatorname*{dom}v\not \supset L_{+}^{2}\supset L_{+}^{\infty}$,
which shows that the inclusion $\operatorname*{dom}v\supseteq L^{2}[0,1]_{+}$,
mentioned in Q7, is not true.

Having in view the texts from Q7, Q8 and Q9, one may wonder what is meant in
\cite{GrOsZa02} by continuity and lower semi\-continuity of $v$ at some point
in $Y$, as well as by local Lipschitzness and subdifferentiability.

In what concerns the local Lipschitzness, it is quite clear that this is meant
in the sense from condition 2.\ in Q4; related to \textquotedblleft
subdifferentiability\textquotedblright, this is not at any point $b$ with
$v(b)\in\mathbb{R}$ as suggested by Q7, but just at a certain point $b$ as in
Q8. As seen in Section \ref{sec2}, the are important differences among the
continuity properties of $g$ and $g|_{\operatorname*{dom}g}$ at points from
$\operatorname*{dom}g$. In fact, inspecting the proof of \cite[Lem.\ 1]%
{GrOsZa02} and the discussion of the modified version of \cite[Examp.\ 1]%
{GrOsZa02}, in \cite{GrOsZa02} one has in view the continuity and the lower
semi\-continuity of $v|_{\operatorname*{dom}v}$ at points in
$\operatorname*{dom}v$.\footnote{Recall that $v|_{\operatorname*{dom}v}$ is
not l.s.c.\ at every $y\in\mathbb{R}\times\mathbb{P}\times\{0\}$ in Example
\ref{ex1}, and $v|_{\operatorname*{dom}v}$ is not l.s.c.\ at every
$y\in\operatorname*{dom}v$ in Example \ref{ex5}.}

Having in view Proposition \ref{p-KZ}, we agree with the remark
\textquotedblleft$v$ is not even continuous there (or anywhere)" from Q8. In
what concerns Q9, on one hand, we would like to see those \textquotedblleft
similar perturbations" which \textquotedblleft show that $v$ is not continuous
anywhere on" $L^{2}[0,1]_{+}\cap\operatorname*{dom}v$; on the other hand, as
already mentioned, we do not agree with the remark \textquotedblleft Of
course, $v$ is lower semicontinuous", which is surely based on \cite[Lem.\ 1]%
{GrOsZa02}.

\section{Some comments on Proposition 2 from \cite{GrOsZa02}}

In Section 6 of \cite{GrOsZa02} one establishes two results on the
Lipschitzness of the value function in infinite-dimensional linear
programming; the second one, Proposition 2, is applied to the assignment model
in \cite[Sect.\ 7]{GrOsZa02}. Our aim is to discuss the proof of
\cite[Prop.\ 2]{GrOsZa02}; for easy reference, we quote its statement and
proof, as well as its preamble:

\medskip

Q10\label{Q10} -- \textquotedblleft Another structural condition is useful for
application to the assignment model. We will use the condition in the context
of a maximization problem and will state it as such.

\smallskip PROPOSITION 2. Let $X$ and $Y$ be \emph{Banach lattices} and let
$A,b$, and $c$ be the data for an LP maximization problem.\footnote{Of course,
it is $c^{\ast}$ instead of $c$.} Assume that

\textbullet\ $A$ is a positive operator which maps the positive cone $X_{+}$
onto $Y_{+}$;

\textbullet\ the order interval $[0,x_{0}]$ is mapped onto the order interval
$[0,Ax_{0}]$ for every $x_{0}\geq0$;

\textbullet\ $A$ is bounded below on the positive cone $X_{+}$, i.e.\ there
exists a constant $M>0$ such that $\left\Vert Ax\right\Vert \geq M\left\Vert
x\right\Vert $ for all $x\geq0$.

Then the value function is Lipschitz on the positive cone $X_{+}$.\footnote{In
fact, it is $Y_{+}$ instead of $X_{+}$.}

\smallskip Proof. Start with $b_{1}\geq0$ and $b_{2}\geq0$. First, consider
the case that $b_{2}\leq b_{1}$. Given $\varepsilon>0$, there is an almost
optimal $x_{1}$ for $b_{1}$, viz.\ \emph{there is} $x_{1}\geq0$ \emph{with}
$Ax_{1}=b_{1}$, and $c^{\ast}(x_{1})+\varepsilon>v(b_{1})$. Since $0\leq
b_{2}\leq b_{1}=Ax_{1}$ and since the positive operator $A$ maps $[0,x_{1}]$
onto $[0,Ax_{1}]$, there is $x_{2}$ such that $0\leq x_{2}\leq x_{1}$ with
$Ax_{2}=b_{2}$. Clearly, $x_{2}$ is feasible for $b_{2}$; hence, $v(b_{2})\geq
c^{\ast}(x_{2})$.

We compute

$v(b_{1})-v(b_{2})\leq c^{\ast}(x_{1})+\varepsilon-c^{\ast}(x_{2}%
)\leq\left\Vert c^{\ast}\right\Vert \left\Vert x_{1}-x_{2}\right\Vert
+\varepsilon\leq\left\Vert c^{\ast}\right\Vert \tfrac{1}{M}\left\Vert
Ax_{1}-Ax_{2}\right\Vert +\varepsilon$

$\ \ \ \ \ \ \ \ \ \ \ \ \ \ \ \ \ \leq\left\Vert c^{\ast}\right\Vert
\tfrac{1}{M}\left\Vert b_{1}-b_{2}\right\Vert +\varepsilon$

\noindent Since this true for arbitrary $\varepsilon>0$, we have that
$v(b_{1})-v(b_{2})\leq\frac{1}{M}\left\Vert c^{\ast}\right\Vert \left\Vert
b_{1}-b_{2}\right\Vert $

\noindent\emph{Switching the roles of} $b_{1}$ \emph{and} $b_{2}$ gives us
$\left\vert v(b_{1})-v(b_{2})\right\vert \leq\frac{1}{M}\left\Vert c^{\ast
}\right\Vert \left\Vert b_{1}-b_{2}\right\Vert $ as desired.

For the general case in which \emph{we do not assume any order dominance
between} $x_{1}$ \emph{and} $x_{2}$, define $x_{3}=x_{1}\wedge x_{2}$. Then
$b_{3}=b_{1}-(b_{1}-b_{2})^{+}$; i.e., $b_{1}-b_{3}=(b_{1}-b_{2})^{+}$. Consequently,

$\left\Vert b_{1}-b_{3}\right\Vert \leq\left\Vert (b_{1}-b_{2})^{+}\right\Vert
\leq\left\Vert b_{1}-b_{2}\right\Vert $.

\noindent Since $0\leq b_{3}\leq b_{2}$ and $v$ \emph{is an increasing
function}, we have that

$v(b_{1})-v(b_{2})\leq v(b_{1})-v(b_{3})\leq c\leq\left\Vert c^{\ast
}\right\Vert \frac{1}{M}\left\Vert b_{1}-b_{2}\right\Vert $

\noindent The same $x_{3}$ works for $v(b_{2})-v(b_{1})$ and we have shown
that $v$ is Lipschitz on $Y$.\footnote{Of course, it must be $Y_{+}$ instead
of $Y$.}" (Our emphasis.)

\medskip

Remarks:

1) Even if not clearly stated, the considered problem is: maximize \ $c^{\ast
}(x)$ s.t.\ $Ax\leq b$ and $x\geq0$; compare with problem (P) on page 273 to
which Proposition 2 is applied. This is also confirmed by the argument
\textquotedblleft Since $0\leq b_{3}\leq b_{2}$ and $v$ \emph{is an increasing
function}, we have that ..." from the end of the proof.

Set $F(b):=\{x\in X\mid x\geq0$, $Ax\leq b\}$ (the feasible set corresponding
to $b\in Y_{+}$).

2) (One had to) Observe first that $F(b)$ is bounded, and so $v(b)\in
\mathbb{R}_{+}$, for every $b\in Y_{+}$

3) By 2) and the definition of $v(b_{1})$, for each $\varepsilon>0$ there
exists $x_{1}\in F(b_{1})$ such that $c^{\ast}(x_{1})+\varepsilon>v(b_{1})$;
hence $x_{1}\geq0$ and $Ax_{1}\leq b_{1}$.

So, why $Ax_{1}=b_{1}$? Without having $Ax_{1}=b_{1}$ one cannot find (using
the hypotheses) $x_{2}\in\lbrack0,x_{1}]$ such that $Ax_{2}=b_{2}$ because
$b_{2}$ could be outside $[0,Ax_{1}]$. How is the argument continued?

4) Assume that for each $\varepsilon>0$ one finds $x_{1}\in F(b_{1})$ such
that $Ax_{1}=b_{1}$ and $c^{\ast}(x_{1})+\varepsilon>v(b_{1})$.
\textquotedblleft Switching the roles of $b_{1}$ and $b_{2}$", will $b_{2}$
have the same property, that is, for each $\varepsilon>0$ one finds $x_{2}\in
F(b_{2})$ such that $Ax_{2}=b_{2}$ and $c^{\ast}(x_{2})+\varepsilon>v(b_{2})$?
\emph{If so, we agree with the estimate }$\left\vert v(b_{1})-v(b_{2}%
)\right\vert \leq\frac{1}{M}\left\Vert c^{\ast}\right\Vert \left\Vert
b_{1}-b_{2}\right\Vert $.

5) 5a) The particular case was the one in which $(0\leq)$ $b_{2}\leq b_{1}$,
that is, the case in which $b_{1}$ and $b_{2}$ are comparable.

5b) Hence, the general case must be \textquotedblleft the one in which we do
not assume any order dominance between" $b_{1}$ and $b_{2}$.

5c) Under 5b), which are $x_{1}$ and $x_{2}$ here? and which is $b_{3}?$ is it
$Ax_{3}$?

5d) We agree with $x_{3}=x_{1}\wedge x_{2}\Rightarrow x_{3}=x_{1}-(x_{1}%
-x_{2})^{+}$. Assume that $b_{k}=Ax_{k}$ for $k\in\{1,2,3\}$ (which could be
envisaged because one had already $b_{k}=Ax_{k}$ for $k\in\{1,2\}$). Because
$b_{3}=b_{1}-(b_{1}-b_{2})^{+}=b_{1}\wedge b_{2}$, one must have
$A(x_{1}\wedge x_{2})=(Ax_{1})\wedge(Ax_{2})$ for $x_{1},x_{2}\in X_{+}$ (or,
equivalently, for $x_{1},x_{2}\in X$). Do the imposed conditions on the data
of \cite[Prop.\ 2]{GrOsZa02} ensure that $A$ is a homomorphism of Banach lattices?

6) Probably, $c$ from the inequality $v(b_{1})-v(b_{3})\leq c$ is $\left\Vert
c^{\ast}\right\Vert \frac{1}{M}\left\Vert b_{1}-b_{3}\right\Vert $, gotten
because $0\leq b_{3}\leq b_{1}$.

\medskip Having in view the above remarks, we consider that the proof of
\cite[Prop.\ 2]{GrOsZa02} needs several clarifications.

So, in our opinion, the authors of \cite{GrOsZa02} did not succeed to
accomplish their goal that emerges from the following text taken from the
beginning of Section 2 of \cite{GrOsZa02}:

\smallskip Q11 -- \textquotedblleft\emph{The present study was motivated by
the problem of showing that there was no gap in the} infinite-dimensional
linear programming \emph{problem that arose in} our studies of \emph{the
continuum assignment problem in [5]}. \emph{The no-gap argument given there
was incomplete; the current paper rectifies that omission}." (Our emphasis.)

\end{document}